\pdfoutput=1 
\documentclass[12pt]{article}
\usepackage{amsmath,amsfonts,amssymb,url,color,epsfig,graphics}
\usepackage[latin1]{inputenc}

\newcommand{\R}{{\mathbb{R}}}

\newcommand{\C}{{\mathbb{C}}}

\newcommand{\N}{{\mathbb{N}}}

\def\ha{\frac{1}{2}}

\def\pa{\partial}
\def\ra{\rightarrow}
\def\ga{\alpha}

\def\ge{\varepsilon}

\def\gl{\lambda}

\def\gs{\sigma}

\def\OPD{pseudo-differential operator}

\newtheorem{coro}{Corollary}[section]
\newtheorem{theo}{Theorem}[section]

\begin{document}

\title{On the spectral theory  of   sign-changing Laplace operators}

\author{Yves  Colin de Verdi\`ere\footnote{Universit\'e Grenoble-Alpes,
Institut Fourier,
 Unit{\'e} mixte
 de recherche CNRS-UGA 5582,
 BP 74, 38402-Saint Martin d'H\`eres Cedex (France);
{\color{blue} {\tt yves.colin-de-verdiere@univ-grenoble-alpes.fr}}}}
\maketitle

\section{Introduction}

The goal of this  note is to describe  a geometric approach to the spectral theory of operators involved in some
meta-materials (see  \cite{C-M}, in particular for references). I became interested in the subject after hearing a lecture by Camille Carvalho
at the ``Journées EDP Rh\^one-Alpes-Auvergne'' last november in Grenoble. Many thanks to the organizers and to Camille. 

The authors of \cite{C-M} study the case of a bounded domain in $\R^d$ with a negatif index of refraction inside
while, in the exterior, the  metric is flat. Then the spectrum is  continuous 
and the interest is about resonances. Using standard black-box scattering methods, as in the book \cite{D-Z}, definition 4.6, 
one can reduce the problem to a compact setting in order to get resonances close to the real eigenvalues of a self-adjoint problem with a discrete spectrum.
This will be our starting point in this
paper. The aim is to show that methods of microlocal analysis can be used in this context in order to get a natural geometric approach.
In particular, we need to use semi-classical Dirichlet-to-Neumann operators for which we were unable to find a reference in the litterature.
This will be the subject of the independent  section \ref{s:dtn}. 
Note that our approach is not limited to the
2D case as in \cite{C-M}. 

\section{Settings and summary  of results}
We denote by 
\begin{enumerate}
\item $X $  a smooth compact  manifold of dimension $d>1$, with a smooth density $|dx|$ and possibly with a boundary $Y$.
We denote by $L^2$ the Hilbert space $L ^2(X,|dx|)$
\item $Z $ a smooth compact
hypersurface in  $X $ with $Z\cap Y=\emptyset $. We assume that  $X \setminus Z=\Omega _+ \cup \Omega _-  $
where the disjoint open sets $\Omega _\pm $ satisfy the following property:
each $z\in Z$ lies in the closures $X_\pm $ of $\Omega _+ $ and of $\Omega _-$. It can be reformulated as
$\pa X_+ \cap (X\setminus Y) = Z $ and $\pa X_- \cap (X\setminus Y )= Z $,
\item  $g_\pm $ some  smooth Riemannian metrics  on  $X_\pm $,
\item  $\Delta _\pm $ the Laplace operators $\Delta _{g_\pm, |dx|}$: if we have   $|dx|=\theta dx_1 dx_2 \cdots dx_d $
in  some local  coordinates $(x_1,\cdots ,x_d)$,
\[ \Delta _\pm =-\theta^{-1}\pa_i \theta g^{ij}_\pm \pa_j ,\]
\item  $k_{+}$ and $k_- $ the restrictions of $g_\pm $ to $Z$.
\end{enumerate}

We consider the operator $M$  defined as follows:
if 
$f=(f_{+},f_- ) $ with $f_\pm $ in  $C^\infty  (X_\pm )$,
\[ M(f_+,f_-)= (\Delta_+ f_+, -\Delta_- f_-) \]
with the transmission   conditions denoted by $({\bf T})$: 
 $f_{+}= f_- $ on $Z$, $f_\pm =0 $ on $Y$ and 
 $\pa _{n_{+}} f_+  = \pa _{n_-} f_- $ on $Z$ where $n_\pm $ are the outgoing unit normals to $Z$ for the metrics $g_\pm $.
 It follows from the Green formula that $M$ is formally symmetric.

 We  denote by $({\bf Ell})$  the  condition 
\begin{equation} \label{equ:ell}
k_{+} <  k_- .\end{equation}
In fact the true condition is that $k_+-k_-$ does not vanish, we will call it $({\bf Ell}^u)$.
 If $Z$ is connected, the sign of $k_+-k_-$ is constant and up to changing $M$
into $-M$, we can assume  $({\bf Ell})$.

We will prove the following results:
\begin{itemize}
\item The condition  $({\bf Ell})$  is  an {\it ellipticity} condition for the transmission condition  $({\bf T})$ on $Z$.
It implies that the spectrum is discrete. As usual we will denote by $\gl_j,~j=1,\cdots , $ the eigenvalues and by $\phi_j,~j=1,\cdots ,$ a
corresponding orthonormal eigenbasis
of eigenfunctions in $L^2(X)$.
\item Weyl formulae for the positive (resp. negative) part of the spectrum.
\item Existence of modes concentrating on $Z$.
\end{itemize}

Note that we use semi-classical analysis with a small parameter $h>0$ as presented for example in the books
\cite{D-S,Zw}. Here we will take $h=h_j=|\gl_j| ^{-\ha }$. In particular, we need to use semi-classical Dirichlet-to-Neumann maps which are  the
subject of the independent section \ref{s:dtn}.

\section{Self-adjointness and ellipticity}

We denote by $g_\pm ^\star $ the dual metrics on the cotangent spaces $T^\star X_\pm $ and by  $k_\pm ^\star $
 the duals of $k_\pm $ defined on $T^\star Z$. 

The operator $M$ is elliptic outside $Y\cup Z$, meaning that the principal symbol $g_+^\star$ (resp. $g_+^\star$)  on $X_+$ (resp. $X_-$) 
 is invertible for all $(x,\xi)\in T^\star (X\setminus (Y\cup Z))\setminus 0$.
Moreover, the Dirichlet boundary condition along $Y$ is well known to be elliptic.
We will show the
\begin{theo}\label{theo:ell}

 If $({\bf Ell}) $ is satisfied, the closure $\bar{M}$ of $M$ has domain
$D\subset  H^2(X_+)\oplus H^2(X_- )$ defined by  the transmission condition $({\bf T })$ and  $\bar{M}$ is self-adjoint. 
The spectrum of $M$ is discrete with eigenfunctions
whose restrictions to  $X_+ $ and to $X_-$ are smooth.

\end{theo}
Note that the previous result is   proved in \cite{C-P-P} in some particular case. Our proof will be different using only
known properties of elliptic boundary values problems.

{\it Proof.}--

The problem is local near $Z$. There are coordinates $(x,y)\in U_+:=]-a,0]\times Z $ (resp. $U_-:=]0,a]\times Z $) with $a>0$ so that 
$g_\pm =dx^2+ K_\pm (x,y,dy ) $ where $K_\pm (x,.,.)$ is a Riemannian metric on $Z$ depending smoothly on $x\leq 0$ (resp. $x\geq 0$).
 The coordinates $x$ are the distances to  $Z$  for $g_\pm$  and the curves
$y=constant$ are geodesics orthogonal to $Z$.

We have  $K_\pm(0,y,dy)=k_\pm(y,dy)$ and
 $n_\pm =\pm \pa _x $.
We will translate the problem in $U=U_+\cup U_-$ into a problem in $U_+$ by using the isomorphism 
$(f_+,f_-)(x,y)\ra (f_+ (x,y), f_-(-x,y))=(\phi_+,\phi_-) $ with the boundary conditions
$\phi_+-\phi_-=0$ on $Z$ and $\pa_x (\phi_+-\phi_-)=0$ on $Z$.
We will check that this mixed Dirichlet-Neumann boundary conditions is elliptic in the sense of \cite{Ho}, definition 20.1.1 .
What H\"ormander says is that we have, for each $(y,\eta)\in T^\star Y\setminus 0$,  to look at the space ${\cal M}$ of solutions of the differential equations
$-\pa_x^2 u_\pm + k_\pm ^\star (y,\eta )u_\pm =0$ which are bounded on $\R^+$ and to check if the  map
$A_{y,\eta} :{\cal M}\ra \C^2$ defined by
$A(u_+,u_-)=(u_+(0)-u_-(0), \pa _x u_+(0)-\pa _x u_- (0))$ is an isomorphism for all $(y,\eta)\in T^\star Y\setminus 0$.
We get $u_\pm(x)=\alpha_\pm  {\rm exp}(-x\sqrt{k_\pm ^\star (y,\eta ) })$   and ellipticity condition is easily seen to be equivalent
to $k_+(y,\eta )\ne k_- (y,\eta )$. We assume  that the difference $k_+ -k_-$ has a constant sign which we can choose  to be
negative without loss of generality.

The operator $M^\star $ which is the adjoint of $M$ has domain $D^\star $ defined as the subset of $L^2$ of pairs $(u_+,u_-)$ so that
$\Delta_\pm u_\pm \in L^2 $ and the condition $({\bf T })$ is satisfied in the sense of distributions. 

Ellipticity of the boundary conditions implies the Fredholm property as proved in \cite{Ho} (see Theorem 20.1.2):
more precisely there exists  a left inverse $A:L^2 \ra D $ and a smoothing  operator $K:D \ra D\cap (C^\infty (X_+)\oplus C^\infty (X_-)) $
so that
$AM^\star = {\rm Id}_{D} + K$.
It follows that $M$ with domain $D$ is closed and self-adjoint and that, for every $\gl $, $\ker (M-\gl ) \subset D\cap (C^\infty(X_+)\oplus C^\infty (X_-)) $.

\section{Weyl formulae}

Let  $N_{+} (\gl )$ (resp. $N_-(\gl)$) be the numbers
$N_{+} (\gl )=\# \{ j|0\leq \gl_j \leq \gl \}$
(resp. $N_{-} (\gl )=\# \{ j|-\gl\leq \gl_j <0 \} $).
let
\[ C_d:={\rm Vol}(B_d)/(2\pi)^d \]
where $B_d$ is the unit ball in $\R^d$. Let us denote by
$V_\pm :={\rm Vol}(X_\pm)$ where the volumes are calculated with the volume forms of $g_\pm $.

Let us remark that the eigenfunctions $\phi_j $ concentrate on $X_{+} $ as $\gl_j \ra +\infty $:
\begin{theo}\label{theo:conc}
if $K\subset X_- \setminus Z$ is a compact set, then, as $\gl_j \ra +\infty $, for all $N\in \N$, 
all derivatives of the $\phi_j $'s are $O\left(\gl_j^{-N}\right)$ on $K$.
\end{theo}
This follows from the fact that $(h^2 \Delta _- +1)\phi_j =0 $ in $X_-$ with $h=\gl_j ^{-\ha}$ and 
$h^2 \Delta _- + 1$ is an elliptic semi-classical operator in $X_-$ (see \cite{Zw}, Chap. 4).

\begin{theo}
We have, as $\gl \ra +\infty $, 
$N_{\pm }(\gl ) \sim C_d V_\pm\gl^{d/2} $. 
\end{theo}

The proof works as follows:
\begin{enumerate}
\item Prove a Weyl asymptotics for $M^2$ which is a positive elliptic operator of degree $4$, namely
\[N_+(\gl)+N_- (\gl)\sim C_d \left(V_++V_-\right)\gl^{d/2}\]
\item Let $X_\ge=X_+\cup B(Z,\ge ) $. Extends  $\phi _j \in L^2_+$ to a function $\psi_j $  with support in $X_\ge $.
Show that the $\psi_j $'s are quasi-modes for $M^2 $ on $X_\ge $.
\item From that, using 1 and minimax,  we get
\[ N^+(\gl) \lessapprox   C_d {\rm Vol}(X_\ge )\gl^{d/2}\]
\end{enumerate}

Let us proceed:

\subsection{Weyl for $M^2$}\label{ss:four}
Let $Q(\phi)=\langle M\phi | M\phi \rangle $ with domain $D_2$ be the quadratic form associated to $M^2$.
We have
\[ H_0^2 (X_+)\oplus H_0^2 (X_-) \subset D_2 \subset H^2 (X_+)\oplus H^2 (X_-)\]
Weyl asymptotics is then well known (see \cite{Ag}, Theorem 4.6): 
\[ \# \{ |\gl_j| \leq \gl \} \sim C_d \left(V_++V_- \right) \gl^{d/2} \]

\subsection{Quasi-modes on $X_\ge $}
From Theorem \ref{theo:conc}, we know that the eigenfunctions $\phi_j \in L^2_+$ are $O(\gl_j ^{-\infty } )$ on every compact $\subset X_- \setminus Z$.
Let us choose $\ge >0$ and $\chi \in C^\infty (\R, [0,1])$ so that $\chi (x)=1$ for $x\leq \ge/2 $ and $\chi (x)=0$
for $x\geq \ge $. Define $\psi_j = \phi_j \chi  $.
We have $(M -\gl_j )\psi_j =O( \gl_j ^{-\infty } )$. Hence the $\psi_j $'s are quasi-modes for $M^2$ on
$X_\ge = X_+\cup B(Z,\ge ) $. 
\subsection{Using minimax}
The dimension of the  space $E_\gl $ generated by the $\psi_j$'s with $\gl_j \leq \gl $ is equivalent to $N_+(\gl )$.
Moreover the quadratic form $Q$ restricted to $E_\gl $ is $\lessapprox \gl  \| .\|_{L^2} $.
Using the minimax, one gets
$N_+ (\gl) \leq C_d \left( V_++ O(\ge ) \right) \gl^{d/2} $.
Moreover we know from Section \ref{ss:four} that
\[ N_+(\gl)+ N_- (\gl ) \sim  C_d\left(V_++V_-\right)\gl ^{d/2} \]
The result follows.

\section{The semi-classical Dirichlet-to-Neumann maps}\label{s:dtn} 

{\it This section is of independent interest. The results are probably known to experts, but we were not able to find a reference.}

Let us consider a compact connected smooth Riemannian manifold $(X,g)$ of dimension $d\geq 2$  with a boundary $Y$.
We give also a smooth density $|dx|$ on $X$.

We will denote by
\begin{itemize}\item $g^\star :T^\star X \ra \R $ the dual of $g$. 
\item 
$k$ the restriction of $g$ to $Y$ and $k^\star:T^\star Y \ra \R $ the dual of $k$
\item  $\Delta =\Delta _{g,|dx|}$  the corresponding Laplace
operator. 
Recall that $\Delta_{g,|dx|}$ is the Friedrichs extension of
the quadratic form $q(f)=\int_X g^\star (x, df(x) ) |dx| $ on $L^2(X,|dx|)$. The principal symbol of  $\Delta_{g,|dx|}$ is $g^\star$.
\item
${\cal E}$ (resp. ${\cal H}$) the elliptic (resp. hyperbolic) region of $T^\star Y$ defined by
\[{\cal E}:= \{ (y,\eta) \in T^\star Y|~ k^\star(y, \eta )>1\} ~,\]
  (resp.
 \[{\cal H}:= \{ (y,\eta) \in T^\star Y|~ k^\star(y, \eta )<1\} ~).\]
 \end{itemize}
We will use semi-classical analysis as presented in several books like \cite{D-S}. Let us recall that $h$ is a small positive  parameter and that we look at
formal expansions in powers of $h$. 

Let us define, for $\ge \in \{ 0, \pm 1 \}$, a semi-classical Dirichlet-to-Neumann map $DtN_\ge$ as an $h$-dependent  operator
on $Y$. If $u$ is a smooth function on $Y$,
one  defines $DtN_\ge(u)$ as follows:
\begin{itemize}
\item If $\ge =0 $ or $\ge=1$,
let $f$ be the unique solution of $(h^2 \Delta  +\ge)f=0, ~f_{|Y}=u $, then
\[ DtN_\ge u= h\frac{\pa f}{\pa n} \]
where $n$ is the outgoing unit normal vector to $Y$.
\item If $\ge =-1$, we will see later that the following definition makes sense:  choose a h-\OPD ~$\Psi $ of degree $0$ on $Y$ so that $\mathrm{WF}_h'(\Psi )$
\footnote{Recall that $WF'(A)$ where $A$ is a pseudo-differential operator is the closure of the complement of the set where the full symbol of $A$
is negligible}
is compactly supported in $ {\cal E}$.
We define then
\[ DtN_{-1}(u) =h\Psi \frac{\pa f}{\pa n} \]
where $f$ is a solution of  $(h^2 \Delta  -1)f=0(h^\infty )$ whose restriction to $Y$ is $\Psi u$.
\end{itemize}

We have then
\begin{theo}\label{theo:dton}
The operators  $DtN_\ge $ are    h-\OPD s on $Y $ whose principal symbol  are $\sqrt{k^\star (y,\eta)+\ge} $ if $\ge=0$ or $\ge=1$, and principal symbol
 $\psi^2  \sqrt{k^\star (y,\eta)-1} $ where $\psi$ is the principal symbol of $\Psi$  if $\ge =-1$.
\end{theo}
{\it  Proof of Theorem \ref{theo:dton}.--}

Let us consider the most difficult case where $\ge=-1$.

{\it Existence.--}  
We will work locally in $Y$: recall that, near each point $y_0\in Y$, there exists local coordinates $(x,y)\in [0,l[ \times U $ with $l>0$ and 
with $U$ an open set of $\R^{d-1}$ so that
$g= dx^2+H(x,y, dy ) $ where $H(x,.)$ is a Riemannian metric on $U$ depending smoothly on $x$.
Let  $u(y)=e^{i\eta y /h}a(y) $ with $a\in C_o^\infty (U)$ and $\eta \in {\cal E}$.
We will construct a WKB solution of
$(h^2 \Delta -1)f= O\left( h^\infty \right)  $ with  
 $f(0,.)=u$.
Here $ O\left( h^\infty \right)  $ means that all derivatives are $O(h^N)$ for all $N$ uniformly in $X$. 

The WKB Ansatz for the function $f$  writes
\[ f(x,y)= e^{iS(x,y) /h} \sum _{j=0}^\infty A_j(x,y) h^j \]
with $S=\eta y + i \mu x + O(x^2) $ where $\mu =\sqrt{ k^\star (y,\eta) -1} $ 
and $\sum _{j=0}^\infty A_j(0,y) h^j=a(y) $.
Recall that
\[ (h^2 \Delta -1)\left( e^{iS(x,y) /h}A(x,y) \right) =  e^{iS(x,y)/h }\left( (g^\star (x,y, d S (x,y))-1) )A + h T(A) + h^2 R(A) \right) \]
with
\begin{itemize}
\item $g^\star $ is the dual of the metric $g$  defined on the cotangent space $T^\star X$. We have
\[ g^\star (x,y,\xi,\eta) = \xi^2 + H^\star (x,y,\eta ) \]
where $H^\star $ is the dual of $H$. The equation $g^\star (x,y, d S (x,y))-1 =0 $ is called the eikonal equation. 
\item $T(A) $ is the transport term: $T$ is a first order differential equation which writes
\[ T(A) = \frac{1}{i} {\cal X}_gA + DA \]
where $D$ is a smooth function, while $ {\cal X}_g $ is the geodesic flow of $g$ restricted to the Lagrangian manifold which is the graph of $\nabla S$:
\[ {\cal X}_g=2 \pa _x S \pa_x +{\cal Y } \]
where ${\cal Y } $ is a vector field depending on $x$ on the manifold $Y$. 
 The transport equation will be
 $T(A)=B$ with $B(0,y)$ given.
\item $R$ is a differential operator of degree $2$.
\end{itemize}

We will construct $S$ and the $A_j$'s as formal series in $x$:
we write $S=\eta_0 y + i\mu x +O(x^2)=  \sum _{j=0}^\infty S_j(y)x^j $
and $A_j=\sum_{j=0}^\infty  A_{j,k} (y) x^k $ with $A_0(0,y)=a(y) $ and $A_j(0,y)=0$ for $j\geq 1$.

 The eikonal equation as well as the transport equation admit unique formal solutions in powers of $x$ because the geodesic flow is formally transverse
 to   the manifold $Z:=\{ (0,y; i\mu, \eta )|y\in U, \mu =\sqrt{ k^\star (y,\eta) -1} \} \subset T^\star X $.

Let us check it for the eikonal equation:
we can write it as
\[ \pa_x S(x,y)  = \Phi (x,y, \nabla _y S (x,y)) \]
with $\Phi $ a smooth complex valued function
$\Phi (x,y,\eta ) =i \sqrt{ H^\star (x,y, \eta ) -1 }$.
Expanding in power series in $x$, we get
\[ \sum _{j=0}^\infty jx^{j-1}S_j(y) = \sum \phi_{k,l\geq 0}(y,\nabla _y S_l) x^{k+l}  \]
We see that the term in $x^{j-1} $ on the righthandside depends only on the $S_l$'s with $l< j$. Hence we can uniquely solve with
$S(0,y)=\eta y $.

We get
\[ (h^2 \Delta -1) (e^{iS(x,y)/h}\sum A_j(x,y)h^j = O(x^\infty ) e^{-\Im S(x,y)/h}=O(h^\infty) \]

Let us denote by $Pu$ the previous extension $f$ of $u(y)=a(y)^{i\eta y/h}$. Note that $Pu$ depends smoothly on $a$ and $\eta$.
Moreover
\[ h\frac{\pa f}{\pa n } _{|x=0} = - h\frac{\pa f}{\pa x } _{|x=0}= e^{iy\eta/h} \left( a(y)\sqrt{ k^\star (y,\eta) -1}
+ \sum_{j=1}^\infty a_j(y) h^j +O(h^\infty) \right) \]  
The righthandside of the previous equation is a $h$-symbol.

We now write the operator $\Psi$ as
\[ \Psi u (y)=\frac{1}{(2\pi h)^{(n-1)} } \int_{U}  \psi (y,\eta)e^{iy\eta/h} \left(\int_U e^{-iy'\eta/h} u(y')dy' \right)  d\eta ~.\]
So that we can write the extension of $\Psi u $ as
\[f= P \Psi u (y)=  \int_{\R^{n-1}} P\left(\psi (y,\eta)e^{iy\eta/h}\right) \hat{u}(\eta ) d\eta \]
where $\hat{u}(\eta )=\frac{1}{(2\pi h)^{(n-1)} }\int _U e^{-i\eta y/h}u(y) dy $ is the $h$-Fourier transform of $u$.
Taking the normal derivative in the integral we get
\[ h\frac{\pa f}{\pa n } _{|x=0}= \int_{\R^{n-1}}\psi (y,\eta)e^{iy\eta/h}\left( \sqrt{ k^\star (y,\eta) -1}+ \sum_{j=1}^\infty a_j(y,\eta) h^j \right) \hat{u}(\eta) d\eta \]
which shows that $DtN u $ as defined by the previous equation is an $h$-\OPD~  of principal symbol 
$\psi^2  \sqrt{ k^\star  -1}$. The extra factor $\Psi$ has no impact on the formula, but is needed in order to get uniqueness in what follows. 

{\it Uniqueness
\footnote{This  argument was suggested to me by Johannes Sj\"ostrand, many thanks to him}.--}

We know, by construction, that   $\mathrm{WF}_h(DtN (u))\subset K \subset {\cal E}$.  Let us give another compact
$K_1$ with $K \subset \dot{K}_1 \subset {\cal E}$. Let us consider then $v=b(y) e^{i\eta' y /h}$ with $\eta '\in K_1 $
and the extension $w $ of $v$ which was build before.
We have
\[ \int _{[0,l[\times U} \left((h^2 \Delta -1)Pf \bar{w} -f(h^2 \Delta -1)\bar{w} \right) |dx| = h^2 \int _U (u \pa_n \bar{v} - \bar{v}\pa_n u ) |dy|  \]
for some smooth density $|dy|$ on $Y$.
Assuming that $u=0$ on $Y$, we get
\[ \int _U  \pa_n u  \bar{v} |dy|  =O(h^\infty ) \]
This holds for any choice of $v$, so that is implies $\mathrm{WF}_h(\pa_n u) \cap  K=\emptyset  $.
Applying the operator $\Psi$, we get
$\Psi\pa_n u=O(h^\infty )$.
This shows that our $DtN$ operator is well defined modulo smoothing operators.

The same kind of calculation shows that $DtN$ is self-adjoint modulo $O(h^\infty )$, so that it can be choosen to be self-adjoint by
modification by a smoothing operator.

\section{Eigenfunctions localized on $Z$ }
This section uses the semi-classical Dirichlet-to-Neumann maps as presented in the  section \ref{s:dtn} which is of independent interest.  

We have  the following result:
\begin{theo}\label{theo:inter}
There exists  a   pseudo-differential equation  $Q_h$ on $\Gamma $ whose principal symbol $q$ is given by 
$q=\ha (k^\star_{+}- k^\star _-) -1  $ so that the solutions of $Q_hu_h=0(h^\infty) $  for small values of $h$   extend to quasi-modes
with eigenvalues $h^{-2}$ of $M$   concentrating on $Z$:
they satisfy $|u_h|=O\left( h^\infty \right)$ on every compact
disjoint of $Z$. 
\end{theo}

{\it Proof.--} 

Let us define $Q_h$ by
\[ Q_h:=\ha (DtN_{+} + DtN_-)(DtN_{+} - DtN_-) \]
whose principal symbol is $\ha(k^\star_{+}- k^\star _- )-1$   
and consider a sequence $(u_j,h_j )$ of solutions of
\[ Q_{h_j} u_j = O(h^\infty ) \]
We have ${\rm WF}(u_j) \subset C:=\{ \|\eta \|_{+}^2 =\|\eta \|_-^2 +2 \} $, so that we are in the elliptic region for $g_{+} $ and we can assume that $\Psi $
is the identity near ${\rm WF}(u_j)$. 
Then, if we define $E=DtN_{+} + DtN_-$ which is elliptic near the  characteristic manifold $C$  of $Q_h$, we get,
by multiplying on the left by a microlocal parametrix  $E^{-1}$ of $E$, 
\[ (DtN_{+} - DtN_-)u_j =O (h^\infty )\]
Extending $u_j$ by $f_j^\pm= P_\pm u_j $, we get solutions
of $(h_j^2 \Delta _\pm \mp 1) f_j^\pm =O(h^\infty ) $.
The restrictions to $Z$ of $f_+$ and $f_-$ coincide while
$\pa_{n^+} f_+ - \pa_{n_-} f_-=r =O(h^\infty )$.
Let us choose $w\in C^\infty (X_+)$ so that $w_{|Z}=0$,
$\pa _{n_+}w= r $ and $w=O(h^\infty )$; we get that
$(f_{+}+w, f_-)$ belongs to $D$ and
\[ (h^2 \Delta _{+}-1) (f_{+}+w)=O(h^\infty),~ (h^2 \Delta _{-}+1) f_{-}=O(h^\infty)\]
So that the pair $(f_{+}+w, f_-)$ is a quasi-mode for $M$.

The pair $(f_j^+,f_j^-)$ is a quasi-mode for the sequence of eigenvalues
$h^{-2}_j $.

\section{Existence of solutions of the equation $Q_hu_h=O(h^\infty )$}

Let us denote by $K$ the Riemannian metric on $Z$ whose adjoint is $(k_+^\star -k_-^\star )/2$ and $\Delta _K $ the corresponding Laplace-Beltrami
operator with an eigenbasis $u_j,~j=1,\cdots $ and corresponding eigenvalues $\mu_j$. Let us define $h_j:=\mu_j^{-\ha}$.   
We want to prove the following result:
\begin{theo} If the basis $(u_j)$ is well choosen,
there exist some sequences of formal power series in $h_j$ of the form $k_j =h_j+\sum_{l=2}^\infty \ga_l h_j^l$ and
$v_j=u_j +\sum_{l=1}^\infty  h_j^l w_{j,l} $ so that
\[ Q_{k_j} v_j= O(k_j^\infty )~.\]
\end{theo}
In particular we have the
\begin{coro}The corresponding quasi-modes of $M$ satisfy the Weyl law associated to $K$.
\end{coro}

Let us prove the previous theorem, assuming first, for simplicity, that the eigenvalues $\gl_j$ of $\Delta _K$ are simple.
Let us look at the first non trivial terms in the expansions, namely $\ga_1$ and $w_{j,1}$. Let us write
$Q_h =h^2 \Delta_K + hP_1 + O(h ^2)$ with $P_1$ a h-\OPD~ of order $0$. 
We get up to terms in $h^2$ the equation
\[ (h_j^2 \Delta_K -1) w_{j,1}=-P_1 v_j -2\ga_1 v_j \]
which can be solved by choosing $\ga_1=-\frac{1}{2}\int _Z v_j P_1 v_j$.
We iterate the previous construction and get the result.

If the some eigenvalues correspponding to $h=h_j $ are degenerate,
we try to get
$v=\sum _{i}x_i u_i + h w $, and similarly we get
\[ (h^2 \Delta -1)w - (2\ga +P_1)\sum _{i}x_i u_i=0 \]
This can be solved iff  the righthandside is orthogonal to all the eigenspace: this says that $-2\ga $ is an eigenvalue
of the matrix $ (<P_1u_i|u_j >) $. This allows to proceed as before.

From the quasi-modes $v_j$ one can built quasi-modes for $M$ by using the extensions of $v_j$ in $X_\pm $.
This ways,  by the usual relation between modes and quasi-modes, one gets exact modes localized on $Z$

Note also that because the $v_j$ are close to the eigenfunctions of $\Delta _K $, we get also quantum ergodicity of these modes
on $Z$ if the geodesic flow of $k$ is ergodic (for the statments and defintions concerning quantum ergodicity, one can look at
\cite{Sh,Zel,CdV}).

\section{Recovering the result of \cite{C-M} for $d=2$}

The case studied in \cite{C-M} is as follows:
$X$  is a bounded open set of $\R^2$ with a boundary $Y$. The metrics 
$g_\pm$ are given by $g_\pm =a_\pm ^{-1} g_0 $ where $g_0$ is the Euclidian metric and $|dx|$ is the  Lebesgue measure coming from $g_0$.
Then $\Delta_\pm = - {\rm div}(a_\pm  \nabla ) $. 
In this case, $Z$ is diffeomorphic to a circle of length $L$ for the metric induced by $g_0$.
We  take the coordinate $y$ as the arc-length on $Z$. The symbol of $Q$ is
\[ \gs (Q)= \frac{a_{+} -a_-}{2} \eta ^2=A(y)\eta^2 \]
It follows that the eigenvalues of $Q$ admits the following asymptotics:
\[ \gl_m \sim \gl_{m+1}\sim \sum_{j=-2}^\infty  c_j m^{-j} \]
with
$ c_{-2}= (2\pi)^2 /l^2 $, 
$l= \int_0^L A^{-\ha} dy $. The parameter $l$  is the length of $Z$ for the metric $A^{-1} dy^2$ associated  to $\gs (Q) $.

\section{Questions}

 What about the interface spectrum  if assumption  $({\bf Ell}^u)$ is not satisfied, but stil the difference of the 2 DtN maps is elliptic of
 one order less? Does essential self-adjointness hold in the later case as it is proved in some cases in \cite{C-P-P}?

In the time dependent case, what  happens on $Z$  to high frequency solutions of the exterior wave equation? How do these waves create some
waves located on $Z$?

\bibliographystyle{plain}

\begin{thebibliography}{}

\end{thebibliography}


\begin{thebibliography}{99}

\bibitem[Ag]{Ag} Shmuel Agmon,
{\it Lectures on Elliptic Boundary Values problems.}
Van Nostrand (1965).



\bibitem[C-M]{C-M}{Camille Carvalho} and {Zo\"is  Moitier},
{\it Scattering resonances in an unbounded transmission problems with sign-changing coefficients.}
IMA J. Appl. Math. {\bf 88 (2)}: 215--257 (2023).


\bibitem[C-P-P]{C-P-P}
Claudio Cacciapuoti, Konstantin Pankrashkin and Andrea  Posilicano, 
{\it Self-adjoint indefinite Laplacians.} 
J. Anal. Math. {\bf 139 (1)}:  155--177 (2019).

\bibitem[CdV]{CdV}
Yves  Colin de Verdi\`ere,
\textit{Ergodicit\'e et fonctions propres du laplacien},
Commun. Math. Phys. {\bf 102} (1985), 497--502.


\bibitem[D-S]{D-S}
Mouez Dimassi and Johannes  Sjöstrand.
{\it Spectral asymptotics in the semi-classical limit.}
London Mathematical Society Lecture Note Series {\bf  268}. 
Cambridge University Press (1999).

\bibitem[D-Z]{D-Z}Semyon Dyatlov and  Maciej Zworski.
Mathematical theory of scattering resonances.
Graduate Studies in Mathematics 200. American Mathematical Society  (2019).

\bibitem[Ho]{Ho} Lars H\"ormander,
{\it The Analysis of Partial Differential Operators III.}
Grundlehren Math. {\bf 274} (2007). 

\bibitem[Sh]{Sh}
Alexander  Shnirelman,
\textit{Ergodic properties of eigenfunctions},
Uspehi Mat. Nauk {\bf 29} (1974), 181--182.

\bibitem[Zel]{Zel}
Steve  Zelditch,
\textit{Uniform distribution of eigenfunctions on compact hyperbolic surfaces},
Duke Math. J. {\bf 55} (1987), 919--941.

\bibitem[Zw]{Zw} Maciej Zworski.
{\it Semi-classical analysis.}
Graduate Studies in Math. {\bf 138}, AMS (2012).


\end{thebibliography}

\end{document}